\documentclass[12pt]{amsart}

\usepackage{amsmath,amssymb,amscd,amsfonts,amsthm,verbatim}
\usepackage[all,cmtip]{xy}
\usepackage{paralist}
\usepackage[mathscr]{eucal}
\usepackage{color}
\usepackage{pdfsync}
\usepackage[]{fontenc}
\usepackage{enumerate}
\usepackage{pgfplots}

\usepackage[pdftex]{hyperref}
\hypersetup{citecolor=blue,linktocpage}

\usepackage[colorinlistoftodos]{todonotes}

\newtheorem{thm}{Theorem}[section]
\newtheorem{lem}[thm]{Lemma}
\newtheorem{prop}[thm]{Proposition}

\newtheorem{cor}[thm]{Corollary}

\newtheorem{assu-nota}[thm]{Assumption--Notation}

\newtheorem*{thm*}{Theorem}

\theoremstyle{definition}
\newtheorem{defn}[thm]{Definition}
\newtheorem{rem}[thm]{Remark}

\newcommand{\R}{\mathbb R}

\numberwithin{equation}{section}

\title{Slope Inequalities for fibrations of non-maximal Albanese dimension}
\author{Miguel \'Angel Barja}

\dedicatory{Dedicated to Professor Fabrizio Catanese for his seventieth birthday}

\address{Miguel \'Angel Barja\\Departament de Matem\`atiques\\Universitat Polit\`ecnica de Catalunya\\ Institut de Matemàtiques de la UPC-BarcelonaTech (IMTech) \\Avda. Diagonal 647\\08028 Barcelona\\Spain}
\email{miguel.angel.barja@upc.edu}

\thanks{The author was supported by MINECO PID2019-103849GB-I00 ``Geometr{\'\i}a Álgebra, Topolog{\'\i}a y aplicaciones multidisciplinares" and by Generalitat de Catalunya SGR2018.}

\begin{document}


\begin{abstract} We study and obtain Slope inequalities for fibred irregular varieties of non-maximal Albanese dimension. We give a comparison theorem between Clifford-Severi and Slope inequalities for this type of fibrations. We also obtain a set of Slope inequalities considering the geometry of the Albanese map and the associated eventual maps.

\end{abstract}

\maketitle
\setcounter{tocdepth}{1}
\tableofcontents
\section{Introduction}

We consider triplets $(X,L,a)$ where $X$ is an irregular, complex, projective variety of dimension $n$, $L$ is a line bundle on $X$ and $a:X\longrightarrow A$ is a non trivial generating map to an abelian variety of dimension $q$. Se will have a fibration $f: X \longrightarrow B$ onto a smooth curve. We will assume that the fibration is {\it irregular}, i.e. ${\rm dim} \, a(F)>0$, where $F$ is a general fibre el $f$.

In this situation, several invariants associated to the triplet can be defined: the {\it continuous rank} $h^0_a(X,L)$, the {\it continuous positive degree} ${\rm deg}_a^+f_*L$ and the {\it eventual map} $\phi_L$. If $h^0_a(X,L)\neq 0$ we define the Clifford-Severi slope of $(X,L,a)$  as $\lambda (L,a)={\rm vol}(L)/h^0_a(X,L)$ and the Slope of $L$ with respect to $f$ as $s(f,L,a)={\rm vol}(L)/{\rm deg}_a^+f_*L$.

When we consider a {\it good} class of triplets ${\mathcal F}$ we denote $\lambda_{\mathcal F}(n)$ and $s_{\mathcal F}(n)$ to be the minimum values of such slopes when $(X,L,a)\in {\mathcal F}$ or $(F,L_{|F},a_{|F})\in {\mathcal F}$, respectively, and $n={\rm dim} \,(X)$.

When varieties in ${\mathcal F}$ are of maximal $a$-dimension, in \cite{B2} we prove that

\smallskip

\begin{equation}\label{main}
{\lambda}_{\mathcal F}(n)\geq \, s_{\mathcal F}(n)\geq \,n\, {\lambda}_{\mathcal F}(n-1).
\end{equation}

\smallskip

We also characterize fibrations with minimal slope in the family. Observe that this provides a way to inductively give higher dimensional Clifford-Severi and Slope inequalities just giving an inequality in low dimension. Moreover, we deduce a huge set of Slope inequalities just using all the existing Clifford-Severi inequalities for varieties of maximal $a$-dimension. In order to obtain these results, in \cite{B2} we develop a version of the method of Xiao adapted to the irregular setting, the so called {\it continuous Xiao's method}.

The aim of this work is to obtain similar results for classes of varieties of non-maximal $a$-dimension. Under this assumption it is easy to see that $\lambda_{\mathcal F}(n)=s_{\mathcal F}(n)=0$, since in this case $h^0_a(X,L)>0$ does not implies bigness. We redefine ${\overline \lambda}_{\mathcal F}(n)$ and ${\overline s}_{\mathcal F}(n)$ when restricting to line bundles $L$ with {\it continuous moving part} $L_c$ big.

Adapting the arguments given in \cite{B2}, our first result is an analogous to (\ref{main}) (see Proposition \ref{pardini}, Theorem \ref{thm1} and Remark \ref{resumen}), and allows to obtain Slope inequalities for varieties of non maximal $a$-dimension from Clifford-Severi inequalities of maximal $a$-dimension ones. More concretely, given a good class of triplets ${\mathcal F}$, we define the subclass ${\mathcal F}_p$ imposing the extra condition that $c(X)={\rm dim}(X)-{\rm dim}\,a(X)\leq p$. Then

\bigskip
\noindent {\bf Theorem A}.
{\it
\begin{itemize}
\item [(i)] ${\overline \lambda}_{{\mathcal F}_p}(n)\geq {\overline s}_{{\mathcal F}_p}(n)\geq (n-p)\,\lambda_{{\mathcal F}_0}(n-p-1).$

\small

\item [(ii)] If equality ${\overline s}_{{\mathcal F}_p}(n)= (n-p)\,\lambda_{{\mathcal F}_0}(n-p-1)$ holds, then $f_*(L\otimes a^*(\alpha))^+$ is semistable for general $\alpha \in {\widehat A}$ and $F$ is covered by $(n-p-1)$-dimensional varieties $V$ of maximal $a_{|V}$-dimension such that $\lambda (R,a_{|V})=\lambda_{{\mathcal F}_0}(n-p-1)$, for some $R\in {\rm Pic} (V)$.
\end{itemize}
}

\bigskip

The technique used here is again the continuous Xiao's method.

\bigskip

The second part of the paper is devoted to obtain new Slope inequalities considering the geometry of  $T$ and $G$, the connected components of the general fibres of $a$ and $\phi_L$, respectively. We refer to Section 2 and Remark \ref{beta} for definitions. Again we use continuous Xiao's method, adapting the arguments of \cite{B} and \cite{J}. Our main result is

\bigskip

\noindent {\bf Theorem B.} {\it Let $f: X \longrightarrow B$ an irregular fibration with general fibre $F$ and $k={\rm dim} \, a(F)$. Then:

\begin{itemize}

\item [(i)] $s(f,L,a)\geq {\rm vol}_{F|G}(L)\, (k+1)!$

\smallskip

\item [(ii)] If $(L_{|F})_c$ is $a_{|F}$-big, then $s(f,L,a)\geq \beta (L_{|T},n,k+1)\,(k+1)!$

\noindent In particular, if $F$ is not uniruled, then $s(f,L,a)\geq 2\,(k+1)!$

\smallskip

\item [(iii)] If equality holds in (i) or (ii), then $f_*(L\otimes a^*(\alpha))^+$ is semistable for general $\alpha \in {\widehat A}$.

\end{itemize}
}

\bigskip

In Section 2 we survey all the techniques we will use, known results on this topic and the involved definitions. Section 3 is devoted to prove theorems A and B.

\bigskip

\noindent  {\underline{Notations and Conventions.}} We work over $\mathbb{C}$. Varieties are projective and smooth unless otherwise stated. We will use notation of divisor or line bundles indistinctly. Given a triplet $(X,L,a)$ we will write $L\otimes \alpha$ instead of $L \otimes a^*(\alpha)$, for $\alpha \in {\widehat A}$.

\bigskip

\noindent {\underline{Acknowledgements.}} The author thanks Lidia Stoppino and Rita Pardini for extremely useful discussions on this topic along the last years.

\bigskip
\bigskip

\section{Preliminaries and technical results}

\bigskip
For benefit of the reader, we collet here a series of preliminaries, definitions, known results and techniques we will use in next section. Main references for these are \cite{B2}, \cite{B}, \cite{BPS3} and \cite{J}.

\smallskip

We consider triplets $(X,L,a)$ where $X$ is a smooth irregular variety of dimension $n$, $a:X \longrightarrow A$ is a nontrivial generating map to an abelian variety of dimension $q$ and $L$ is a line bundle on $X$ such that $h^0_a(X,L)={\rm min}\{h^0(X,L\otimes \alpha)\,|\,\alpha \in {\widehat A}\}\neq 0$.

\bigskip

{\noindent {\underline {\it Multiplication maps}}}. We will offen consider situations {\it up to a multiplication map}, meaning that we will consider base changes via a multiplication map on $A$ by some $d$, which is \'etale of degree $d^{2q}$:

\begin{equation}\label{multiplicationmap}
\xymatrix{
X^{(d)}\ar[d]_{a_d}\ar[r]_{\nu_d} &X\ar[d]^a\\
A\ar[r]_{\nu_d}&A}
\end{equation}

We will denote $L^{(d)}:=\nu_d^*(L)$. Continuous rank and volume are multiplicative through a multiplication map.

\bigskip

\noindent {\underline {\it Continuous moving divisor $L_c$}}.  Up to a blow-up, there is a decomposition $L=P+W$ such that, for any $d>>0$ and divisible and any general $\alpha \in {\widetilde A}$, $P^{(d)}$ is base point free and is the moving divisor of $|L^{(d)}\otimes \alpha|$ and $W^{(d)}$ is its fixed divisor. Following \cite{B}, $P$ and $W$ are called the {\it continuous moving divisor} and {\it continuous fixed divisor} of $L$, respectively. According to the notation of \cite{J}, we will set $L_c:=P$ for the continuous moving part.

\bigskip

\noindent {\underline {\it Eventual map and eventual degree.}} Up to a blow-up, there is a factorization of the map $a$, $X \rightarrow X_L\rightarrow A$ such that the map $\phi_L: X\longrightarrow X_L$ verifies the following properties:

\begin{itemize}
\item $L_c=\phi_L^*(R_L)$ for some line bundle $R_L$ on $X_L$ which induces a base point free generically finite morphism on $X_L$ (\cite{BPS3}).
\item Up to a multiplication map, the linear system $|L_c^{(d)}\otimes \alpha|$, for $\alpha$ general, is base point free and induces the map $\phi_L^{(d)}: X^{(d)}\longrightarrow X_L^{(d)}$ (\cite{B}, \cite{BPS3}).
\item Since the map $\phi_L$ factorizes $a$, it is generically finite provided $X$ is of maximal $a$-dimension. It is birational if ${\rm deg}\, a=1$.
\item The map $\phi_L$ is generically finite provided $L_c$ is $a$-big (\cite{B}).
\item A birational model of the map $\phi_L$ is given by the natural map $\rho: X \dashrightarrow \mathbb{P}_A(a_*L)$, where $X_L:=\rho (X)$ (\cite{J}).
\item When the eventual map $\phi_L$ is generically finite, we define the {\it eventual degree} $L$ to be $e(L)={\rm deg} (\phi_L)$ (cf. \cite{BPS3}, Section 3). We can extend the definition to any $L$ just considering the degree of the finite part in the Stein factorization of $\phi_L$.
\item $\kappa (L_c)={\rm dim}\phi_L(X) \geq {\rm dim} \, a(X)$ (\cite{J}).
\end{itemize}

\bigskip

\noindent {\underline {\it Good classes of triplets.}} Given a family ${\mathcal F}$ of triplets with $h^0_a(X,L)\neq 0$, we say that the family is {\it good} if it stable via the following four operations:

\begin{itemize}

\item[(1)] If $(X,L,a)\in {\mathcal F}$, then $({\overline X},\sigma ^*L,a\circ \sigma)\in {\mathcal F}$, where $\sigma:{\overline X}\longrightarrow X$ is a birational morphism.
\item[(2)] If $(X,L,a)\in {\mathcal F}$, then $(X^{(d)},L^{(d)},a_d)\in {\mathcal F}$.
\item[(3)] If $(X,L,a)\in {\mathcal F}$, then $(X,L',a)\in {\mathcal F}$ for $L'\leq L$ such that $h^0_a(X,L')>0$.
\item[(4)] If $(X,L,a)\in {\mathcal F}$, then $(M,L_{|M},a_{|M})\in {\mathcal F}$, for a general smooth $M$, moving in a base point free linear system on $X$.
\end{itemize}

\bigskip

\noindent {\underline {\it Clifford-Severi inequalities (maximal $a$-dimension)}}. Given a triplet $(X,L,a)$ with $h^0_a(X,L)>0$, we define its Clifford-Severi slope as

$$\lambda (L,a)=\frac{{\rm vol}(L)}{h^0_a(X,L)}$$

\noindent which remains constant under multiplication maps. Given a good class ${\mathcal F}$ we define ${\lambda}_{\mathcal F}(n)$ to be the infimum of the Clifford-Severi slopes of triplets in ${\mathcal F}$, of dimension $n$. Clifford Severi-Inequalities for a given good class ${\mathcal F}$ are inequalities of type

$$
{\lambda}_{\mathcal F}(n)\geq \lambda_n.
$$

\medskip

There are many known Clifford-Severi inequalities for different good classes of maximal $a$-dimension varieties: see Remark 4.7 in \cite{B2} for a (almost) complete list. For example:

\begin{itemize}

\item Higher dimensional Severi inequality states that ${\lambda}_{\mathcal F}(n)\geq 2\,n!$ if ${\mathcal F}$ is defined by the property that $L$ is {\it numerically subcanonical}.

\item For a general $L$ we have ${\lambda}_{\mathcal F}(n)\geq e(L)\,n!$ (\cite{BPS2}).

\end{itemize}

\bigskip

\noindent {\underline {\it Clifford-Severi inequalities (non-maximal $a$-dimension)}}. In the case of triplets of non-maximal $a$-dimension, situation is not so clear and depends heavily on conditions of bigness of $L$ or $L_c$ and the geometry of the fibre of the map $a$ or $\phi_L$. The case of irregular threefolds is well understood by results of Zhang (\cite{Z3}).

Here you can find a (non complete) list of known results for arbitrary dimension $n$. We set $k={\rm dim} a(X)$, and $G$ and $T$ for a connected component of the general fibre of $\phi_L$ and $a$, respectively.

In \cite{B}, Main Theorem and Remark 5.8, the author proves

\begin{itemize}
\item If $L_c$ is $a$-big and $L$ is numerically $r$-subcanonical, then $\lambda(L,a)\geq \delta (r)\,k!$.
\item If $L$ is nef and $a$-big then $\lambda(L,a)\geq (L_{|G})^{n-k}\,k!\geq \, k!$.
\end{itemize}

\noindent When $k=n-1$, Zhang gives a better bound (\cite{Z2}):

\begin{itemize}
\item If $g$ is the genus of the curve $T$, then $\lambda (L,a)\geq 2\, \frac{g-1}{g+n-2}\, n!.$
\end{itemize}

\noindent Finally Jiang (\cite{J}) gives a set of inequalities depending on the geometry of $G$ or $T$. The simplest ones are the following (see Proposition 3.6 and Theorem 3.1 in \cite{J} and Remark \ref{beta} for a more detailed result):

\begin{itemize}
\item If $L$ is big, then $\lambda(L,a)\geq {\rm vol}_{X|G}(L)\,k!$.
\item If $L_c$ is big and $T$ is not uniruled, then $\lambda(L,a)\geq 2 \, k!$.
\end{itemize}

\begin{rem} The proof of Main Theorem (iii) in \cite{B} uses implicitly the volume of $L_{|G}$ (see Remark 5.8 in the cited paper). In Corollary B (iii), loc. cit., extending this inequality to $K_X$ in the singular setting, it is erroneously assumed that ${\rm vol}_{X|G}(K_X)\geq 1$, assuming that the minimal variety $X$ is Gorenstein.

\end{rem}

\bigskip

\noindent {\underline {\it Irregularly fibred triplets.}} Given a triplet $(X,L,a)$, we will say that it is {\it irregularly fibred} if moreover we have a fibration $f:X\longrightarrow B$ onto a smooth curve, such that ${\rm dim}\,a(F)>0$, where $F$ is a general smooth fibre of $f$.

If $f$ is an irregular fibration as above, the family of vector bundles $\{f_*(L\otimes \alpha)\}_{{\alpha}\in {\widehat A}}$ has constant type of Harder-Narashimann filtration for $\alpha \in U_0$, for some nonempty open set $U_0$. We set $\{(r_i,\mu_i)\}$ for theit ranks and slopes.

For $\alpha \in U_0$ we will write $f_*(L\otimes \alpha)^+$ for the biggest nef subbundle of $f_*(L\otimes \alpha)$ and we denote

$${\rm deg}_a^+f_*L={\rm deg}f_*(L\otimes \alpha)^+\geq 0.$$

We will define the (continuous) Slope of $L$ w.r.t. $f$ to be:

$$
s(f,L,a)=\frac{{\rm vol}( L)}{{\rm deg}_a^+f_*L}\in (0+\infty].
$$

\noindent which is also constant under multiplication maps.

\bigskip

\noindent {\underline {\it Slope inequalities.}} Given a good class ${\mathcal F}$, we will say that an irregular fibration $f$, with general fibre $F$, is of type ${\mathcal F}$ if $(F,L_{|F},a_{|F})\in {\mathcal F}$ (the triplet $(X,L,a)$ is not necessarily in ${\mathcal F}$). We define $s_{\mathcal F}(n)$ to be the infimum of slopes of fibrations $f$ of type ${\mathcal F}$, where $n={\rm dim} X$.

Slope inequalities for the family ${\mathcal F}$ is a set of inequalities for any $n$:

$$
s_{\mathcal F}(n) \geq \lambda_n.
$$

\medskip

In \cite{HZ} Hu and Zhang give slope inequalities for $L=K_f$ and $X$ of maximal Albanese dimension, by direct computation, giving properties of the limit cases.

In \cite{B2}, Theorem 4.11, a broad generalization is given, for any $L$, establishing an equivalence between Clifford-Severi inequalities and Slope inequalities for a given good class ${\mathcal F}$ of varieties of maximal $a$-dimension. Moreover, this result allows to produce automatically a whole set of Clifford-Severi and Slope inequalities for any dimension, just given one inequality in low dimension, typically 1 or 2 (see Remark 4.12 and Corollary 1.1 in \cite{B2}). The main result can be stated as:

$$
\lambda_{\mathcal F}(n)\geq s_{\mathcal F}(n)\geq n\,\lambda_{\mathcal F}(n-1).
$$

\bigskip

\noindent {\underline {\it Continuous Xiao's method and derived inequalities.}} Take a general $\alpha_0\in {\widehat A}$ and let $L_0=L\otimes \alpha_0$. After a suitable blow-up and multiplication map there is a filtration by nef line bundles:

$$T_1 \leq T_2 \leq ...\leq T_m  \leq L_0$$

\noindent such that, for all $i$ $N_i:=T_i-\mu_iF$ is nef and, if $P_i:={N_i}_{|F}$, then

\begin{itemize}
\item $P_i$ is a base point free linear system on $F$ such that $h^0(F,P_i)=h^0_{a_{|F}}(F,P_i)\geq r_i$.
\smallskip
\item ${\rm deg}_a^+f_*(L)={\rm deg}({\mathcal E}_m^{\alpha_0})=\sum_{i=1}^{m}r_i(\mu_i-\mu_{i+1})$
\end{itemize}

\noindent By convention we can take, coherently, $(N_{m+1},\mu_{m+1})=(T_m,0)$ or, in case $L$ is nef, $(N_{m+1},\mu_{m+1})=(L,0)$.

Fix  $r\leq n$. Consider an ordered, increasing partition of the set $\{ \, 1,...,m \, \}$ given by subsets $I_s$, $s=1,...,r-1$ (some of the sets $I_s$ may be empty), with $I_{r-1}\neq \emptyset$. Define, decreasingly, for $s=1,...,r-1$

$$b_s=\{
           \begin{array}{cc}
             {\rm min} I_s & {\rm if} \,\, I_s \neq \emptyset \\
             b_{s+1} & {\rm otherwise} \\
           \end{array}$$

\bigskip

 Then we have that, for any $Q_1,...,Q_{n-r}$ nef $\mathbb{Q}$-Cartier divisors the following inequality holds:

\begin{equation}\label{Xiaogeneral}
Q_1...Q_{n-r}\left[N^r_{m+1}-(\sum _{s=1}^{r-1} (\prod_{k>s}P_{b_k})\sum_{i \in I_s}(\sum_{l=0}^{s}P_i^{s-l}P_{i+1}^{l})(\mu_i-\mu_{i+1}))\right]\geq 0.
\end{equation}

\bigskip

We will use the following particular cases.

Taking  $r=n$, and any partition:

\begin{equation}\label{xiaobuena}
{\rm vol}(L)={\rm vol} (L_0)\geq N_{m+1}^n\geq \sum _{s=1}^{n-1} (\prod_{k>s}P_{b_k})\sum_{i \in I_s}(\sum_{l=0}^{s}P_i^{s-l}P_{i+1}^{l})(\mu_i-\mu_{i+1}).
\end{equation}

Taking any $r$, the trivial partition $I_{r-1}=I$ ($I_s=\emptyset$ for $s<r-1$) and $Q_i=N_{m+1}$:

\begin{equation}\label{xiaoalbanesenomaxima}
{\rm vol}(L)={\rm vol} (L_0)\geq N_{m+1}^n\geq 
\sum_{i=1}^m P_m^{n-r} \left[P_{i+1}^{r-1}+P_{i+1}^{r-2}P_i+...+P_i^{r-1}\right](\mu_i-\mu_{i+1}).
\end{equation}

\bigskip

\section{Slope inequalities for non maximal $a$-dimension fibrations}

\medskip

In \cite{B2} we study the equivalence of Slope and Clifford-Severi inequalities for general irregular varieties and fibrations in a good class ${\mathcal F}$, which mostly applies to maximal $a$-dimension varieties. Our aim is to study whether this equivalence can be extended to classes of varieties of non-maximal $a$-dimension.

To this aim, we need to make this condition stable by smooth sections (condition (3)), and the right condition to impose is on the codimension $c(X)=\dim X-\dim a(X)$, which we called condition $Q_p$, i.e. triplets such that $c(X)\leq p$.

Parts of the main result in \cite{B2} given in Theorem 4.11 hold for families of non maximal $a$-dimension as well, but they are not interesting since in these cases $\lambda_{\mathcal F}(n)$ and $s_{\mathcal F}(n)$ vanish. The reason is that condition $h^0_a(X,L)>0$ does not imply bigness if the variety is not of maximal $a$-dimension. Indeed, if $a:X\longrightarrow A$ verifies that $k=\dim a(X)<n=\dim X$, take $L=a^*H$ for any $H$ very ample on $A$. Then clearly ${\rm vol}(L)=0$ and $h^0_a(X,L)>0$.The same phenomena occur if $X$ is fibred: we can construct examples with $\deg_a^+f_*L\neq 0$ and ${\rm vol}(L)=0$.
So we need to impose extra hypotheses to obtain nontrivial inequalities. Natural conditions are $a$-bigness of $L$ or its continuous moving part $L_c$.  There are several strategies (see \cite{B} and \cite{J}), according to whether $L$ or $L_c$ are $a$-big. Observe that bigness of $L_c$ implies bigness of $L$ but the viceversa does not hold (see Remarks 3.7 and 3.8 in \cite{B}). Observe also that in \cite{B} it is shown that $a$-bigness of $L_c$ implies bigness of $L$. In particular, for $L_c$, bigness and $a$-bigness are equivalent, provided $h^0_a(X,L)\neq 0$.

If we restrict our good families ${\mathcal F}$ adding the condition of bigness of $L$ (or $L_c$), the resulting subfamily is not {\it good} since bigness is not stable by subbundles (and so condition (4) fails).

Nevertheless, we can obtain some closely related Slope inequalities for fibrations of non maximal $a$-dimension with adapted arguments.

Let us first introduce some extra notation.

\begin{defn} Given a class ${\mathcal F}$ of triplets we define:

\begin{itemize}
\item [(i)] ${\mathcal F}_p=\{(X,L,a)\in {\mathcal F}\,|\, c(X)\leq p\,\}.$
\smallskip
\item[(ii)] ${\overline \lambda}_{{\mathcal F}_p}(n)={\rm inf}\{\lambda(L,a)\,|\, (X,L,a)\in {\mathcal F}_p,\,n=\dim X,\,\, L_c \, {\rm big}\}.$
\smallskip
\item[(ii')] ${\widehat \lambda}_{{\mathcal F}_p}(n)={\rm inf}\{\lambda(L,a)\,|\, (X,L,a)\in {\mathcal F}_p,\,n=\dim X,\,\, L \, {\rm big}\}.$
\smallskip
\item[(iii)] ${\overline s}_{{\mathcal F}_p}(n)={\rm inf}\{s(f,L,a)\,|\, f \,{\rm of}\,\,{\rm type}\,\, {\mathcal F}_p,\,n=\dim X, \,\,(L_{|F})_c \,{\rm big}\}.$
\smallskip
\item[(iii')] ${\widehat s}_{{\mathcal F}_p}(n)={\rm inf}\{s(f,L,a)\,|\, f \,{\rm of}\,\,{\rm type}\,\, {\mathcal F}_p,\,n=\dim X, \,\,(L_{|F}) \,{\rm big}\}.$

\end{itemize}

\end{defn}

\medskip
\begin{rem} \begin{itemize}
\item If the class ${\mathcal F}$ is good, so is ${\mathcal F}_p$.
\item If we consider classes of maximal $a$-dimension as in \cite{B2}, then ${\mathcal F}={\mathcal F}_0$. In this case, since $h^0_a(X,L)\neq 0$ implies bigness, we have that $\lambda_{\mathcal F}(n)={\overline {\lambda}}_{{\mathcal F}_0}(n)={\widehat {\lambda}}_{{\mathcal F}_0}(n)$ and $s_{\mathcal F}(n)={\overline s}_{{\mathcal F}_0}(n)={\widehat s}_{{\mathcal F}_0}(n)$.
\end{itemize}

\end{rem}

\bigskip

One of the two inequalities between Clifford-Severi and Slope inequalities given in Theorem 4.11 in \cite{B2} holds without change in this new setting:

\begin{prop}\label{pardini} Let ${\mathcal F}$ be a good class of triplets of irregular varieties. Then, for all $n$ and $p\leq n-2$ we have

$$
{\overline \lambda}_{{\mathcal F}_p}(n)\geq {\overline s}_{{\mathcal F}_p}(n)\,\,\,\, {\rm and} \,\,\,\, {\widehat \lambda}_{{\mathcal F}_p}(n)\geq {\widehat s}_{{\mathcal F}_p}(n).
$$
\end{prop}

\begin{proof}
We refer to the proof of the maximal $a$-dimension case using Pardini's trick given in Theorem 4.11, (ii) of \cite{B2}. As pointed out in Remark 4.13 (loc.cit.), only properties (1), (2) and (3) of a good class are used, and bigness of $L$ or $L_c$ are maintained in all the process.

The condition $p\leq n-2$ ensures that ${\rm dim} \, a(X)\geq 2$. In this case, the sections of $\nu_d^*(H)$ are irreducible and $f_d$ has connected fibres (and so it is an irregular fibration).

The rest of the proof holds without changes.




\end{proof}

\bigskip

The reverse inequality is more subtle and does depend heavily on bigness properties of $L_c$ or $L$. When $L_c$ is $a$-big two possible strategies are possible. The first option, following \cite{B}, is by hyperplane section argument since the eventual map allows us to maintain the process inside the good class ${\mathcal F}$.  This is the content of Theorem \ref{thm1} and Remark \ref{resumen}.

The second option is to consider the geometry of $T$, the connected component of the general fibre of $a_{|F}$. This approach uses Theorem 3.1 in \cite{J}, adapted to the irregularly fibred case via a suitable use of the continuous Xiao's method. This is the content of Theorem \ref{thm2} (ii).

In general $(L_{|F})_c$ may not be $a_{|F}$-big. In this case, we can also obtain a good lower estimation of $s(f,L,a)$, but we need to consider the geometry of $G$, a connected component of  the general fibre of the eventual map $\phi_{L_{|F}}$. This approach adapts the argument in Proposition 3.6 in \cite{J} to the relative setting via continuous Xiao's method and is the content of Theorem \ref{thm2} (i).



Observe that bounds in Theorem \ref{thm2} in general are sharper than those in Theorem \ref{thm1} when considering a single fibred triplet $(X,L,a)$ but strongly depend on properties not well behaved in a good class ${\mathcal F}$.

As a by product of the use of Xiao's method, in theorems \ref{thm1} and \ref{thm2} we can give properties in the limit cases, being those in Theorem \ref{thm1} analogous to those obtained in the cases of maximal $a$-dimension varieties.
\bigskip


\begin{thm}\label{thm1}

Let ${\mathcal F}$ be a good class of triplets. Let $(X,L,a)$ be a fibred triplet $f: X\longrightarrow B$ of type ${\mathcal F}$ such that $f$ is of $a$-dimension $k$ (i.e., $k={\rm dim}\,a(F)$). Assume that $L_c$ is $f$-big. Then:

\begin{itemize}

\item [(i)] $s(f,L,a) \geq \,(k+1)\,{\lambda}_{{\mathcal F}_0}(k).$

\smallskip

\item[(ii)] If equality holds, then:
   \begin{itemize}
   \item There is a family of varieties $V$ of dimension $k$ covering $F$, and line bundles $R\leq L_{|V}$ such that $(V,R,a_{|V})$ are of maximal $a_{|V}$-dimension and verify the Clifford-Severi equality ${\rm vol}(R)={\lambda}_{{\mathcal F}_0}(k)\,h^0_{a_{|V}}(V,R).$
   \item $f_*(L\otimes \alpha)^+$ is semistable for general $\alpha \in {\widehat A}$. If, moreover, $L$ is nef, then $f_*(L\otimes \alpha)$ is nef and semistable.
   \end{itemize}
\end{itemize}
\end{thm}

\bigskip

\begin{rem}\label{thm1simple}
Statement (i) of the above theorem, combined with Clifford-Severi inequalities for maximal $a$-dimension varieties as given in Remark 4.9 in \cite{B2} gives a broad set of Slope inequalities.

For example, in general we have $s(f,L,a) \geq (k+1)!$ and, if $L$ is numerically $r$-subcanonical, then $s(f,L,a) \geq \delta(r)(k+1)!$

\end{rem}

\medskip

\begin{rem} \label{resumen}
Theorem \ref{thm1} together with Proposition \ref{pardini} can be rephrased as, for $p\leq n-2$:

    $$
    {\overline \lambda}_{{\mathcal F}_p}(n)\geq {\overline s}_{{\mathcal F}_p}(n)\geq (n-p)\,\lambda_{{\mathcal F}_0}(n-p-1).
    $$

\end{rem}

\bigskip

\begin{proof}

\noindent (i) Let $(X,L,a)$ be a triplet with a fibration $f: X \longrightarrow B$ of type ${\mathcal F}$.

Let us apply inequality (\ref{xiaoalbanesenomaxima}) with $r=k+1={\rm dim}\,a(F)+1$:

\begin{equation} \label{uno}
{\rm vol}(L)={\rm vol} (L_0)\geq N_{m+1}^n\geq \sum_{i=1}^m P_m^{n-k-1}\left[P_{i+1}^{k}+P_{i+1}^{k-1}P_i+...+P_i^{k}\right](\mu_i-\mu_{i+1}).
\end{equation}

Since $L_c \leq T_m=N_{m+1}$ (see Lemma 3.1 in \cite{B2}), we have by hypothesis that $P_m={T_m}_{|F}$ is big and hence $ a_{|F}$-big. Hence, its eventual map is generically finite. Moreover, since it is continuously globally generated by construction, up to a multiplication map we can assume that the linear system $|P_m|$ is base point free. Take $V_1,...,V_{n-k-1}\in |P_m|$ general sections and let $V=V_1\cap ... \cap V_{n-k-1}$. Hence $V$ is a smooth variety of maximal $a_{|V}$-dimension $k$. Let $R_i={P_i}_{|V}$. By the properties of a good class, we have that the triplet $(V,R_i,a_{|V})\in {\mathcal F}_0$. Hence we have that

\begin{equation}\label{dos}
R_i^{k}\geq {\lambda}_{{\mathcal F}_0}(k)\,h^0_{a_{|V}}(V,R_i).
\end{equation}

Observe that

$$
P_m^{n-k-1}P_i^{k}=R_i^{k}.
$$

Finally, we use that $P_i-P_m\leq 0$ and then $h^0_{a_{|F}} (F,P_i-P_m)=0$, and the same holds by cutting by successive $V_i$. Then we can conclude that

\begin{equation}\label{tres}
h^0_{a_{|V}}(V,R_i)\geq h^0_{a_{|F}}(F,P_i)\geq r_i.
\end{equation}

\noindent Finally, observe that using general Clifford-Severi inequality for irregular varieties (Main Theorem in \cite{B})

\begin{equation}\label{cuatro}
\delta_i:=R_{i+1}^{k}+R_{i+1}^{k-1}R_i+...+R_i^{k}\geq R_{i+1}^{k}+k\,R_i^{k}\geq \lambda_{{\mathcal F}_0}(k)\,( r_{i+1}+kr_i)\geq (k+1)\,\lambda_{{\mathcal F}_0}(k)\,r_i.
\end{equation}

Then we can conclude

\begin{equation} \label{cinco}
{\rm vol}(L)={\rm vol} (L_0)\geq N_{m+1}^n\geq \sum_{i=1}^m (k+1)\,\lambda_{{\mathcal F}_0}(k)\,r_i(\mu_i-\mu_{i+1})=(k+1)\,\lambda_{{\mathcal F}_0}(k)\,{\rm deg}_a^+f_*L.
\end{equation}



\medskip

\noindent (ii) Assume that equality holds. Then we have equality in (\ref{dos}), (\ref{tres}), (\ref{cuatro}) and (\ref{cinco}), which imply, for all $i=1,...,m$:

\begin{itemize}
\item $r_{i+1}=r_i$,
\item $h^0_{a_{|V}}(V,R_i)= h^0_{a_{|F}}(F,P_i)= r_i$,
\item $R_i^{k}={\lambda}_{\mathcal F}(k)_0\,h^0_{a_{|V}}(V,R_i)$,
\item ${\rm vol} (L_0)=T_m^{n}$.
\end{itemize}

Hence we have that $m=1$ (so $f_*(L\otimes \alpha)^+$ is semistable), $h^0_{a_{|V}}(V,R_m)= h^0_{a_{|F}}(F,P_m)= r_m$, and $(V,R_m)$ verifies the Clifford-Severi equality $R_m^{k}={\lambda}_{{\mathcal F}_0}(k)\,h^0_{a_{|V}}(V,R_m)$.

Observe that if equality holds then $L$ is big since $s(f,L,a)>0$. We have that $L_0=T_m+Z_m$  and ${\rm vol} (L_0)=T_m^{n}$. If $L$ is moreover nef, then we have that $Z_m=0$ (cf. Theorem A in \cite{FKL}). Hence:

$${\rm rank}f_*(L_0)^+=r_1=h^0_{a_{|F}}(F,P_1)=h^0_{a_{|F}}(F,{L_0}_{|F})={\rm rank}f_*(L_0)$$

\noindent and hence $f_*(L_0)$ is semistable and nef.
\end{proof}

\bigskip

\begin{lem}\label{fujita}  Let $X$ be a smooth, projective variety and $L$ a big line bundle on $X$. Let $\phi: X \longrightarrow Y$ be a fibred space, $G$ a general fibre of $\phi$ and $R\in {\rm Pic} (Y)$ a line bundle such that $L'=\phi^* (R)\leq L$. Then:

$$
{\rm vol}_X (L)\geq {\rm vol}_{X|G}(L)\,{\rm vol}_{Y}(R).
$$
\end{lem}

\medskip

\begin{proof}
We have that $L=\phi ^*(R)+Z$ with $Z\geq 0$. The result is obvious if $L$ is nef. We will reduce to this case via Fujita Approximation theorem. In the general case, assume $R$ is big, otherwise the result is trivial. Following Theorem 3.5 in \cite{BDPP}, there is an extension of the volume function given by the moving intersection numbers which is non decreasing in each factor, superadditive and coincides with the intersection product for nef line bundles. Let $e={\rm dim} (Y)$. For any birational compatible modifications of $X$ and $Y$ and any decompositions  $L=W_1+E_1$ and $R=W_2+E_2$ such that $W_i$ are big and nef $\mathbb{Q}$-divisors and $E_i$ are effective $\mathbb{Q}$-divisors, we have:

$$
{\rm vol}(L)=<L,....,L>\,\geq\, <W_1,...,W_1,\phi^*(W_2),...,\phi^*(W_2)>=W_1^{n-e}(\phi^*(W_2))^e=({W_1}_{|G})^{n-e}W_2^e
$$

\noindent since the moving intersection numbers are nondecreasing and $W_i$ are nef.

To conclude, we apply Fujita Approximation theorem for the volume and the relative volume (see, for example, Proposition 2.11 and Theorem 2.13 in \cite{elmnp}). For any $\epsilon >0$, there are birational modifications of $X$ and $Y$ (that we can make compatible with the above hypotheses) and decompositions $L=W_1+E_1$ and $R=W_2+E_2$ as above and such that

\begin{itemize}

\item ${\rm vol}_{X|G}(L) \geq {W_1}_{|G}^{n-e}\geq {\rm vol}_{X|G}(L)-\epsilon$ and
\item ${\rm vol}_Y(R)\geq W_2^e \geq {\rm vol}_Y(R)-\epsilon$.
\end{itemize}

We apply the above inequality for any  such decompositions and we conclude that

$$
{\rm vol}_X (L) \geq {\rm vol}_{X|G}(L)\, {\rm vol}_Y (R).
$$
\end{proof}

\bigskip

\begin{rem}\label{beta}
Let $Z$ be a smooth, projective variety, and $L$ a line bundle on $Z$.  In \cite{J}, two invariants to compute the positivity of $L$ are defined: $\delta(L)$ and $\delta_1(L)$, the second being the minimum of volumes of subline bundles of $L$ inducing generically finite maps, when restricted to general positive dimensional subvarieties $V$ covering $Z$. We clearly have that $\delta_1(L)\geq cov.gon(Y)$. For any $k\leq s\leq {\rm dim} (Z)$, define

$$
\beta(L,s,k)={\rm min}\{\binom{s}{k}\delta (L),\delta_1(L)\}.
$$

\noindent  We have that $\beta(L,s,k)\geq 1$ and that $\beta(L,s,k)\geq 2$, provided $Y$ is not uniruled.

In \cite{J} Theorem 3.1, the following Clifford-Severi inequality is proved. Consider a triplet $(X,L,a)$ of dimension $n$. Let $k={\rm dim} \, a(X)$ and $T$ be a connected component of the general fibre of $a$. Assume that $L_c$ is $a$-big. Then:

$$
\lambda(L,a)\geq \beta(L_{|T},n,k)\,k!
$$

\end{rem}

\bigskip

\bigskip

\begin{thm}\label{thm2}
Let $(X,L,a)$ be a fibred triplet. Let $G$ be a connected component of the general fibre of the eventual map $\phi_{L_{|F}}$ and let $T$ be a connected component of the general fibre of $a_{|F}$. Then

\begin{itemize}

\item [(i)] If $L_{|F}$ is $a_{|F}$-big, then $s(f,L,a)\geq {\rm vol}_{F|G}(L)\,(k+1)!$

\smallskip

\item [(ii)] If $(L_{|F})_c$ is $a_{|F}$-big, then $s(f,L,a)\geq \beta (L_{|T},n,k+1)\,(k+1)!$

\noindent In particular, if $F$ is not uniruled, then $s(f,L,a)\geq 2\,(k+1)!$

\smallskip
\item [(iii)] If equality holds in (i) or (ii), then $f_*(L\otimes \alpha)^+$ is semistable for general $\alpha \in {\widehat A}$.

\end{itemize}
\end{thm}

\bigskip

\begin{proof}  (i) Since $(T_m)_{|F}=P_m$ is continuously globally generated and coincides by construction with $(L_{|F})_c$, up to a multiplication map and a birational modification, we can consider that $|P_m|$ induces the eventual map of $L_{|F}$ and $a_{|F}$ factorizes through this map. Up to a further birational modification and multiplication map, we can consider the relative map induced by the quotient $f_*(L_0)^+=f_*(T_m)^+={\mathcal E}_m\longrightarrow T_m$. The Stein factorization of such map gives a relative fibration $\phi_m: X\longrightarrow X_m$ over $B$, with general fibre $G$ and a factorization as $a={\overline a}\circ \phi_m$ for some ${\overline a}: X_m\longrightarrow A$. We can assume $X_m$ to be smooth. Then $T_m=\phi_m^*(L_m)$ for some line bundle $L_m$ on $X_m$ such that $g_*(L_m)^+=f_*(T_m)^+=f_*(L_0)^+={\mathcal E}^{\alpha_0}_m$, where $g: X_m \longrightarrow B$ is the induced fibration over $B$.

 When restricted to a general fibre ${\overline F}$, this is just the Stein factorization of the eventual map $\phi_{L_{|F}}$.







Then we apply Lemma \ref{fujita}:

$$
{\rm vol}_X (L) \geq {\rm vol}_{X|G}(L)\, {\rm vol}_{X_m}(L_m).
$$

Observe that ${L_m}$ restricted to the fibre of $g$ is ${\overline a}$-big. Hence, we can apply Remark \ref{thm1simple} to $(X_m,L_m,{\overline a})$ and obtain

$$
{\rm vol}_{X_m}(L_m)\geq (k+1)!\,{\rm deg}_{\overline a}^+g_*(L_m)=(k+1)!\,{\rm deg}_a^+ f_*(L_0).
$$

\bigskip

\noindent (ii) Consider the set of indexes $I=\{1,....,m\}$ and let us construct an increasing ordered partition as follows. For $s=1,...,n-1$, consider $I_s=\{i\in I \,|\, \kappa (P_i)=s\,\}$. Recall that $\kappa (P_i)={\rm dim} \, \phi_{P_i}(F)$. Since $(L_{|F})_c=P_m$ is $a_{|F}$-big then its eventual map is generically finite, and so $I_{n-1}\neq \emptyset$. On the other hand, since $a_{|F}$ factorizes through any eventual map, we have that ${\rm dim}\,\phi_i(F)\geq k$, and hence $I_s=\emptyset$ if $s<k$.

Consider now Xiao's inequality (\ref{Xiaogeneral}):

\begin{equation}\label{xiaobuena}
N^n_{m+1}\geq \sum _{s=1}^{n-1} (\prod_{k>s}P_{b_k})\sum_{i \in I_s}(\sum_{l=0}^{s}P_i^{s-l}P_{i+1}^{l})(\mu_i-\mu_{i+1}).
\end{equation}

\bigskip

Consider the Stein factorization of the eventual maps induced by $|P_i|$, $\phi_i: F \longrightarrow F_i$, and let $R_i\in {\rm Pic}(F_i)$ be such that $P_i=\phi_i^*(R_i)$. We can assume that all the $F_i$ are smooth and that $\phi_i$ factorizes through $\phi_j$ if $i<j$. Let $G_i$ the generic fibre of $\phi_i$. Observe that if $i,i'\in I_s$, then $G_i=G_{i'}$ and so we denote it by $G_s$. Observe that ${\rm dim} \, G_s=n-1-s.$ Since eventual maps factorizes $a_{|F}$, we also have maps $a_i: F_i\longrightarrow A$ such that $a_{|F}=a_i\circ \phi_i$.

Then we have that, for $i\in I_s$:

$$
\sum_{l=0}^{s}P_i^{s-l}P_{i+1}^{l}\geq (s+1)P_i^s\geq \left[(k+1)R_i^s\right]G_s.
$$

Since $(F_i,R_i,a_i)$ is a triplet such that $(R_i)_c=R_i$ is big and of $a_i$-dimension $k$, we can apply general Clifford-Severi inequaliy to obtain $R_i^s\geq k! h^0_{a_i}(F_i,R_i)=k!h^0_{a_{|F}}(F,P_i)\geq k! \, r_i$.

Summing up we have

\begin{equation}\label{paraigualdad}
\sum_{l=0}^{s}P_i^{s-l}P_{i+1}^{l}\geq \left[k! \,(r_{i+1}+kr_i)\right]\,G_s\geq \left[(k+1)!\,r_i\right]\,G_s
\end{equation}



When $s\leq n-2$, we have that $(\prod_{k>s}P_{b_k})G_s$ is the volume of $P_{b_{n-1}}$ restricted to a curve in $T$. Since  $|P_{b_{n-1}}|$ induces a generically finite map, we have that $(\prod_{k>s}P_{b_k})G_s\geq \delta_1(L_{|T})$. Hence, for any $i\in I_s$ we have

$$
(\prod_{k>s}P_{b_k})(\sum_{l=0}^{s}P_i^{s-l}P_{i+1}^{l})(\mu_i-\mu_{i+1}))\geq (s+1)\delta_1(L_{|T})k!r_i \geq (k+1)!\delta_1(L_{|T})r_i.
$$

When $s=n-1$, and $i\in I_{n-1}$ we can use Theorem 3.1 in \cite{J} to obtain

$$
P_i^{n-1-l}P_{i+1}^l\geq P_i^{n-1}\geq \beta(L_{|T},n-1,k) k! r_i,
$$

\noindent and so

$$
(\sum_{l=0}^{n-1}P_i^{n-1-l}P_{i+1}^{l})(\mu_i-\mu_{i+1})\geq n\beta(L_{|T},n-1,k)\,k!\,r_i.
$$

Now we take the minimal lower bound for $s\leq n-2$ and $s=n-1$ and observe that

$$
{\rm min}\{n\,\beta (L_{|T},n-1,k)\,k!, (k+1)!\,\delta_1(L_{|T})=\beta (L_{|T},n,k+1)(k+1)!
$$

We finally obtain

$$
{\rm vol}(L)\geq N^n_{m+1}\geq  \beta(L_{|T},n,k+1) (k+1)! \sum_{i=1}^mr_i(\mu_i-\mu_{i+1})=\beta(L_{|T},n,k+1)(k+1)!\,{\rm deg}_a^+f_*L.
$$

\bigskip

\noindent (iii) If equality holds in (i), then equality holds for $s(g,L_m,{\overline a})$ and so Theorem \ref{thm1} applies.
If equality holds in (ii), then equality holds in any step and in particular in formula (\ref{paraigualdad}). If $r_{i+1}=r_i$ for all $i$ then $m=1$ and $f_*(L_0)^+$ is semistable.
\end{proof}

\bigskip

     \end{document}